\documentclass[bezier,11pt]{article}
\usepackage{amssymb}
\usepackage{latexsym}
\usepackage{tikz}

\newtheorem{thm}{Theorem}
\newtheorem{ob}[thm]{Observation}
\newtheorem{prop}[thm]{Proposition}
\newtheorem{lem}[thm]{Lemma}

\newtheorem{conj}{Conjecture}
\newtheorem{cor}[thm]{Corollary}

\newcommand{\oc}{{\rm oc}}
\newcommand{\barX}{{\overline{X}}}

\newcommand{\cN}{{\cal N}}
\newcommand{\cC}{{\cal C}}
\newcommand{\cD}{{\cal D}}

\newcommand{\cub}{{\rm cubic}}

\newcommand{\ind}{{\rm ind}}

\newcommand{\smallqed}{{\tiny ($\Box$)}}
\newcommand{\qed}{$\Box$}

\newcommand{\proof}{\noindent\textbf{Proof. }}

\def\vertex(#1){\put(#1){\circle*{4}}}
\def\vertexo(#1){\put(#1){\circle{2}}}
\def\vert(#1){\put(#1){\circle*{1.5}}}
\def\verto(#1){\put(#1){\circle{1.5}}}
\def\lab(#1)#2{\put(#1){\makebox(0,0)[c]{#2}}}
\newenvironment{unnumbered}[1]{\trivlist
\item [\hskip \labelsep {\bf #1}]\ignorespaces\it}{\endtrivlist}

\begin{document}
\bibliographystyle{plain}

\title{Induced Cycles in Graphs}
\author{$^1$Michael A. Henning, $^2$Felix Joos, \\
$^2$Christian L\"{o}wenstein, and $^2$Thomas Sasse
 \\
\\
$^1$Department of Mathematics \\
University of Johannesburg \\
Auckland Park, 2006 South Africa\\
E-mail: mahenning@uj.ac.za \\
\\
$^2$Institute of Optimization and Operations Research,\\
Ulm University,\\
Ulm 89081, Germany,\\
E-mail: felix.joos@uni-ulm.de,\\
E-mail: christian.loewenstein@uni-ulm.de,\\
E-mail: thomas.sasse@uni-ulm.de}

\date{}
\maketitle

\begin{abstract}

The maximum cardinality of an induced $2$-regular subgraph of a graph $G$ is denoted  by $c_\ind(G)$.
We prove that if $G$ is an $r$-regular graph of order~$n$, then $c_\ind(G) \geq \frac{n}{2(r-1)} + \frac{1}{(r-1)(r-2)}$
and we prove that if $G$ is a cubic claw-free graph on order~$n$, then $c_\ind(G) > 13n/20$ and this bound is asymptotically best possible.
\end{abstract}

{\small \textbf{Keywords:} Induced regular subgraph; induced cycle; claw-free; matching; $1$-extendability }\\
\indent {\small \textbf{AMS subject classification: 05C38, 05C69}}

\section{Introduction}

The problem of finding a largest induced $r$-regular subgraph of a given graph for any value of $r \ge 0$ has attracted much interest and dates back to Erd\H{o}s, Fajtlowicz, and Staton~\cite{er}. Cardoso et al.~\cite{cakalo} showed that it is NP-hard to find a maximum induced $r$-regular subgraph of a given graph. Lozin et al.~\cite{lomopu} established efficient algorithms for special graph classes including $2P_3$-free graphs, while Moser and Thilikos~\cite{moth} studied FPT-algorithms for finding regular induced subgraphs.

The special case of finding a largest induced $r$-regular subgraph when $r=0$ is the  well-studied problem of finding a maximum independent set in a graph.
When $r = 1$, the problem is to find a maximum induced matching in a graph which has also received considerable attention in the literature. In this paper, we focus our attention on finding a largest induced $2$-regular subgraph of a given graph.
The situation turned out to be much more complex than the independent set problem and the induced matching problem due to the fact that the components of a $2$-regular subgraph are no longer bounded.  Local arguments and techniques which earlier proved successful do not suffice and global arguments need to be found that take into account structural properties of the graph.

The \emph{induced cycle number} $c_\ind(G)$ of a graph $G$ is the maximum cardinality of an induced $2$-regular subgraph of $G$.
Our aim in this paper is threefold.

Our first aim is to establish NP-hardness of $c_{\rm ind}(G)$ for graphs $G$ of maximum degree~$4$.
Our second aim is to provide a lower bound on the induced cycle number of a general graph in terms of its order, size and maximum degree. As a consequence of this result, we obtain a lower bound on the induced cycle number of a regular graph.
Our third aim is to establish an asymptotically best possible lower bound on the induced cycle number of a cubic claw-free graph.
Our proof techniques rely heavily on matching results and intricate counting arguments.

\subsection{Notation}

Let $G$ be a graph with vertex set $V(G)$ and edge set $E(G)$, and of order~$n(G) = |V(G)|$ and size $m(G) = |E(G)|$.
Let $v$ be a vertex in $V(G)$.
The \emph{open neighborhood} of $v$ is $N_G(v) = \{u \in V \, | \, uv \in E(G)\}$ and the \emph{closed neighborhood of $v$} is $N_G[v] = \{v\} \cup N(v)$. The degree of $v$ is $d_G(v) = |N_G(v)|$.
If $d_G(v) = k$ for every vertex $v \in V(G)$, we say that $G$ is a $k$-regular graph.
A $3$-regular graph is also called a \emph{cubic graph}.
Let $\Delta(G)$ be the maximum degree of $G$.

A \emph{path} on $n$ vertices is denoted by $P_n$.
A complete graph $K_3$ we call a \emph{triangle}.
For a subset $S \subseteq V$, the subgraph induced by $S$ is denoted by $G[S]$.
If $X$ and $Y$ are two vertex disjoint subsets of $V$, then we denote the set of all edges of $G$ that join a vertex of $X$ and a vertex of $Y$ by $E(X,Y)$.
The complete graph on four vertices minus one edge is called a \emph{diamond}.
We say that a graph is \emph{$F$-free} if it does not contain $F$ as an induced subgraph. In particular, if $F = K_{1,3}$, then we say that the graph is \emph{claw-free}. An excellent survey of claw-free graphs has been written by Flandrin, Faudree, and Ryj\'{a}\v{c}ek~\cite{ffr}.

Two edges in a graph $G$ are \emph{independent} if they are vertex disjoint in $G$.
A set of pairwise independent edges of $G$ is called a \emph{matching} in $G$.
If every vertex of $G$ is incident with an edge of a matching $M$, then $M$ is a \emph{perfect matching} in $G$.
A graph $G$ is $1$-\emph{extendable} if given an arbitrary edge~$e$ in $G$ there exists a perfect matching in $G$ that contains $e$.
Matchings in graphs are extensively studied in the literature (see, for example, the book by Lov\'{a}sz and Plummer~\cite{lopl} and the survey articles by
Plummer~\cite{Pl03} and Pulleyblank~\cite{Pu95}).

A \emph{block} of a graph $G$ is a maximal $2$-connected subgraph of $G$.
A block $B$ is a \emph{trivial block} if $B=K_2$
and $B$ is a \emph{cycle block} if $B$ is a cycle.
A \emph{cactus} is a graph containing trivial blocks and cycle blocks only.

\section{Complexity Results}

The basic complexity question concerning the decision problem for $2$-regular subgraphs takes the following form:

\medskip
\textbf{INDUCED CYCLE}: \\
\indent  \textbf{Instance:} A given graph $G$ and a positive integer $k$. \\
\indent  \textbf{Question:} Is $c_\ind(G)\ge k$?

\medskip
The \emph{independence number} $\alpha(G)$ of a graph $G$ is the maximum cardinality of an independent set in $G$. The problem of finding a maximum independent is set is called the maximum independent set problem and is a classic NP-hard optimization problem.

\medskip
\textbf{INDEPENDENT SET PROBLEM}: \\
\indent  \textbf{Instance:} A given graph $G$ and a positive integer $k$. \\
\indent  \textbf{Question:} Is $\alpha(G) \ge k$?

\medskip
The INDEPENDENT SET PROBLEM is NP-complete for the class of cubic planar graphs~\cite{GaJo77}.

\begin{thm}
INDUCED CYCLE is NP-complete for planar graphs of maximum degree at most~$4$.
\end{thm}

\proof
It is easy to see that INDUCED CYCLE is in NP.
In order to show that INDUCED CYCLE is NP-complete, we reduce the INDEPENDENT SET PROBLEM to INDUCED CYCLE.
Assume we want to decide whether the independence number of a cubic planar graph $G$ is at least $k$.
For every vertex $v$ of $G$, let $e_v$ be an arbitrary edge incident with $v$.
Let $G'$ arise from $G$ by subdividing every edge $e_v$ once and if $e_v=e_u$ for two distinct vertices $u,v$, then subdividing $e_v$ twice.
For $v\in V(G)$, let $v'$ be the vertex of degree $2$ that is adjacent to $v$ and arises from subdividing the edge $e_v$.
Let $H$ arise from $G'$ and $n(G)$ many distinct paths $\{P_v\}_{v\in V(G)}$ of order $3n(G)-2$
by joining for every $v\in V(G)$ one endvertex of $P_v$ with $v$ and the other endvertex with $v'$.
Let $C_v$ be the induced cycle containing $P_v$ in $H$.
Note that $H$ is planar, has maximum degree $4$, and $n(H)=3(n(G))^2$.

If $Z$ is a $2$-regular subgraph of $H$,
then a vertex in $P_v$ is contained in $Z$ if and only if all vertices of $P_v$ are contained in $Z$ if and only if $C_v$ is a component of $Z$.
Let $Z$ be a $2$-regular subgraph of $H$ such that the order of $Z$ is maximum.
Assume for contradiction that $Z$ contains a component $C\not= C_v$ for every $v\in V(G)$.
Thus $C$ contains at most two vertices  of $C_v$ for $v\in V(G)$.
Let $v$ be such that $C_v\cap C \not=\emptyset$.
Replacing $C$ by $C_v$ results a $2$-regular subgraph of $H$ of larger order, which is a contradiction to our assumption that the order of $Z$ is maximum.

Therefore, we may assume that $Z= \bigcup_{v\in S}C_v$.
On the one hand, if $uv\in E(G)$, then by construction $C_v$ and $C_u$ cannot be both in $Z$.
Thus $S$ is an independent set.
On the other hand, for every independent set $I$, the graph $ \bigcup_{v\in I}C_v$ is a $2$-regular subgraph of $G$.
Hence, there is an independent set $I$ of order $k$ in $G$ if and only if there is a $2$-regular subgraph of order $3kn(G)$ in $H$.~\qed

\section{General Bounds}

In this section, we provide a lower bound on the induced cycle number of a general graph in terms of its order, size and maximum degree.

\begin{thm}
If $G$ is a graph, then
\[
c_\ind(G)\ge \frac{m(G)-n(G)+1}{(\Delta(G)-2)(\Delta(G)-1)}.
\]
 \label{t:general_bd}
\end{thm}
\proof
Consider the following Greedy-Algorithm.
Remove an induced cycle together with all its neighbors from $G$ until the resulting graph has no cycles any more.
Let $F$ denote the graph that arises by applying the Greedy-Algorithm on $G$
and let $C_1,\ldots,C_t$ denote the cycles of $G$ that are chosen by the Greedy-Algorithm.
Then, $\bigcup_{i=1}^tV(C_i)$ is an induced $2$-regular subgraph in $G$.
We set $\ell=\sum_{i=1}^t|V(C_i)|$, which implies that $c_\ind(G)\ge \ell$.
Let $N$ be the number of vertices that are removed from $G$ by the Greedy-Algorithm
but do not belong to any of the cycles $C_1,\ldots,C_t$.
Then the following equation and inequalities holds.
\begin{eqnarray}
n(F)&=&n(G)-\ell-N\label{eq:1}\\
m(F)&\ge& m(G)-\ell-N\cdot\Delta(G)\label{eq:2}\\
N&\le& \ell\cdot(\Delta(G)-2)\label{eq:3}.
\end{eqnarray}
Since $F$ is a forest, we conclude that
\begin{eqnarray}
m(F)&\le &n(F)-1\label{eq:4}.
\end{eqnarray}
Therefore,
\begin{eqnarray*}
\ell(\Delta(G)-2)(\Delta(G)-1)&
\stackrel{(\ref{eq:3})}{\ge}&N\cdot\Delta(G)-N\\
&\stackrel{(\ref{eq:2}),(\ref{eq:1})}{\ge}&(m(G)-m(F)-\ell)+(n(F)-n(G)+\ell)\\
&\stackrel{(\ref{eq:4})}{\ge}& m(G)-n(G)+1.
\end{eqnarray*}
The desired result now follows from our earlier observation that $c_\ind(G)\ge \ell$.~\qed

\bigskip
As an immediate consequence of Theorem~\ref{t:general_bd} we obtain the following lower bound on the induced cycle number of a regular graph.

\begin{cor}\label{d regular}
For $k \ge 3$, if $G$ is a $k$-regular graph,
then
$$c_\ind(G)\ge \frac{n(G)}{2(k-1)}+\frac{1}{(k-2)(k-1)}.$$
\end{cor}

\begin{cor}
If $G$ is a cubic graph of order~$n$, then $c_\ind(G) > n/4$.
 \label{c:cubic_bd}
\end{cor}
Note that for $k\geq 4$, we have $c_\ind(K_{k,k})=\frac{2}{k}$ and the cubic graph $G$ depicted in Figure~\ref{f:cubic} has induced cycle number $\frac{n(G)}{2}$.
For every $k$, we do not find any graph where $\frac{c_\ind(G)}{n(G)}$ is smaller.

\section{Claw-Free Cubic Graphs}

In this section we establish an asymptotic tight lower bound on the induced cycle number of a connected claw-free cubic graph. We shall prove the following result.

\begin{thm}
If $G$ is a cubic claw-free graph of order~$n$, then $c_\ind(G) > 13n/20$.
 \label{t:clawfree}
\end{thm}

\subsection{Preliminary Results and Known Results}

A component of a graph is \emph{odd} or \emph{even} depending on whether its order is odd or even, respectively. 
For a graph $G$, let $\oc(G)$ denote the number of odd components of $G$.
We shall need the following well-known matching result due to Tutte~\cite{Tutte}.

\begin{thm}{\rm (Tutte's Theorem)}
A graph $G$ has a perfect matching if and only if $\oc(G-S) \le |S|$ for every proper subset $S \subseteq V(G)$.
\end{thm}

\noindent
Moreover, we use the following result.

\begin{lem}{\rm (\cite{HeLo12})}
The vertex set of a connected cubic claw-free graph $G \ne K_4$ can be uniquely partitioned into sets each of which induce a triangle or a diamond in $G$.
 \label{l:triangle_diamond}
\end{lem}

\begin{figure}[htb]
\tikzstyle{every node}=[circle, draw,fill=black!0, inner sep=0pt,minimum width=.175cm]
\begin{center}
\begin{tikzpicture}[thick,scale=.475]
  \draw(0,0) { 
    +(3.64,8.30) -- +(2.66,7.89)
    +(1.17,6.41) -- +(0.76,5.43)
    +(0.76,3.32) -- +(1.17,2.34)
    +(2.66,0.86) -- +(3.64,0.45)
    +(5.74,0.45) -- +(6.72,0.86)
    +(8.21,2.34) -- +(8.61,3.32)
    +(8.61,5.43) -- +(8.21,6.41)
    +(6.72,7.89) -- +(5.74,8.30)
    +(4.69,8.75) -- +(4.69,7.81)
    +(3.64,8.30) -- +(4.69,8.75)
    +(4.69,8.75) -- +(5.74,8.30)
    +(5.74,8.30) -- +(4.69,7.81)
    +(4.69,7.81) -- +(3.64,8.30)
    +(2.66,7.89) -- +(1.59,7.47)
    +(1.59,7.47) -- +(1.17,6.41)
    +(1.17,6.41) -- +(2.26,6.81)
    +(2.26,6.81) -- +(1.59,7.47)
    +(2.66,7.89) -- +(2.26,6.81)
    +(0.76,3.32) -- +(1.25,4.38)
    +(1.25,4.38) -- +(0.00,4.38)
    +(0.00,4.38) -- +(0.76,3.32)
    +(0.00,4.38) -- +(0.76,5.43)
    +(0.76,5.43) -- +(1.25,4.38)
    +(6.72,7.89) -- +(7.12,6.81)
    +(6.72,7.89) -- +(7.78,7.47)
    +(7.78,7.47) -- +(7.12,6.81)
    +(7.78,7.47) -- +(8.21,6.41)
    +(8.21,6.41) -- +(7.12,6.81)
    +(8.61,5.43) -- +(8.13,4.38)
    +(8.13,4.38) -- +(9.38,4.38)
    +(9.38,4.38) -- +(8.61,5.43)
    +(8.13,4.38) -- +(8.61,3.32)
    +(8.61,3.32) -- +(9.38,4.38)
    +(7.12,1.94) -- +(8.21,2.34)
    +(8.21,2.34) -- +(7.78,1.28)
    +(7.78,1.28) -- +(7.12,1.94)
    +(7.12,1.94) -- +(6.72,0.86)
    +(6.72,0.86) -- +(7.78,1.28)
    +(5.74,0.45) -- +(4.69,0.94)
    +(4.69,0.94) -- +(3.64,0.45)
    +(3.64,0.45) -- +(4.69,0.00)
    +(4.69,0.00) -- +(4.69,0.94)
    +(4.69,0.00) -- +(5.74,0.45)
    +(2.66,0.86) -- +(2.26,1.94)
    +(2.26,1.94) -- +(1.17,2.34)
    +(1.17,2.34) -- +(1.37,1.06)
    +(1.37,1.06) -- +(2.26,1.94)
    +(1.37,1.06) -- +(2.66,0.86)
    +(0.76,3.32) node{}
    +(1.25,4.38) node{}
    +(1.17,2.34) node{}
    +(4.69,0.94) node{}
    +(4.69,0.00) node{}
    +(7.12,1.94) node{}
    +(7.78,1.28) node{}
    +(7.12,6.81) node{}
    +(7.78,7.47) node{}
    +(0.76,5.43) node{}
    +(2.26,6.81) node{}
    +(1.59,7.47) node{}
    +(4.69,7.81) node{}
    +(4.69,8.75) node{}
    +(8.13,4.38) node{}
    +(9.38,4.38) node{}
    +(2.26,1.94) node{}
    +(0.00,4.38) node{}
    +(1.37,1.06) node{}
    +(2.66,0.86) node{}
    +(3.64,0.45) node{}
    +(5.74,0.45) node{}
    +(6.72,0.86) node{}
    +(8.21,2.34) node{}
    +(8.61,3.32) node{}
    +(8.61,5.43) node{}
    +(8.21,6.41) node{}
    +(6.72,7.89) node{}
    +(5.74,8.30) node{}
    +(3.64,8.30) node{}
    +(2.66,7.89) node{}
    +(1.17,6.41) node{}
  };
 \end{tikzpicture}
\end{center}
\vskip -0.6 cm \caption{A diamond-necklace $N_8$} \label{f:N8}
\end{figure}
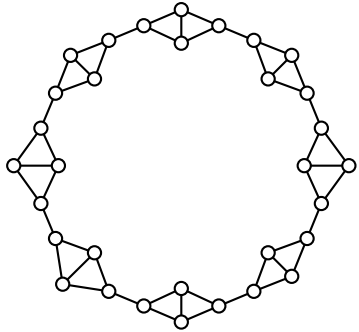

For $k \ge 2$ an integer, let $N_k$ be the connected cubic graph
constructed as follows. Take $k$ disjoint copies $D_1, D_2, \ldots,
D_{k}$ of a diamond, where $V(D_i) = \{a_i,b_i,c_i,d_i\}$ and where
$a_ib_i$ is the missing edge in $D_i$. Let $N_k$ be obtained from the
disjoint union of these $k$ diamonds by adding the edges
$\{a_ib_{i+1} \mid i = 1,2,\ldots,k-1\}$  and adding the edge
$a_{k}b_1$. Following the notation in~\cite{HeLo12}, we call $N_k$ a \emph{diamond-necklace with $k$
diamonds}. Let $\cN_\cub = \{ N_k \mid k \ge 2\}$. A
diamond-necklace, $N_8$, with eight diamonds is illustrated in
Figure~\ref{f:N8}.

\begin{ob}
If $G \in \cN_\cub$ is a diamond-necklace of order~$n$, then $c_\ind(G) = 3n/4$.
 \label{necklace}
\end{ob}

We shall need the following result.

\begin{thm}
Every connected cubic bridgeless graph is $1$-extendable.
 \label{t:matching}
\end{thm}
\proof
Let $G$ be a connected cubic bridgeless graph and let $e=uv\in E(G)$ be an arbitrary edge of $G$.
We will show that $G'=G-\{u,v\}$ has a perfect matching. Since $G$ does not contain parallel edges, we note that $\delta(G')\ge 1$. Among all proper subset of vertices of $G'$, let $X \subset V(G')$ be chosen so that

\medskip
\hspace*{0.5cm} (1) $\oc(G'-X) - |X|$ is maximized, \\
\indent \hspace*{0.5cm} (2) Subject to (1), $|X|$ is maximized.

\medskip
\noindent
We proceed further with a series of claims.

\begin{unnumbered}{Claim~1}
Every component in $G'-X$ is an odd component.
\end{unnumbered}
\proof Suppose to the contrary that $G$ has an even component. Let $x$ be a vertex of such an even component and let $X' = X \cup \{x\}$. Then, $\oc(G'-X') \ge \oc(G'-X) + 1$ and $|X'| = |X| + 1$, implying that $\oc(G'-X') - |X'| \ge \oc(G'-X) - |X|$, contradicting our choice of the set $X$.~\smallqed


\begin{unnumbered}{Claim~2}
$|E(X,\barX)| \ge 3\oc(G'-X)-4$.
\end{unnumbered}
\proof Let $s=\oc(G'-X)$ and let $G_1,\ldots,G_s$ denote the components of $G'-X$.
Let $\barX = V(G') \setminus X$.
For $1\le i\le s$, let $V_i = V(G_i)$.
By Claim~1, every component $G_i$ is odd, and so $\barX = \bigcup_{i=1}^s V_i$.
For $1\le i\le s$, we define $t_i=3|V_i|-\sum_{v \in V_i} d_{G'}(v)$.
Since four edges join $uv$ and the vertices of $G'$, we conclude $\sum_{i=1}^s t_i \le 4$.
By using the fact that there are no edges joining $V_i$ and $V_j$ for $i\not=j$, we obtain for $1\le i\le s$
\[
|E(X,V_i)| = \sum_{v \in V_i} d_{G'}(v) - \sum_{v \in V_i} d_{G_i}(v) =
3|V_i|-t_i-2|E(G_i)|.
\]
Since $|V_i|$ is odd, we note that $|E(X,V_i)|$ and $t_i$ have different parity.

\begin{unnumbered}{Subclaim~2.1}
$|E(X,V_i)| \ge 3-t_i$ for all $1\le i\le s$.
\end{unnumbered}
\proof Suppose $t_i=0$ and $|E(X,V_i)|=1$ for some $i$, $1\le i\le s$. Then the unique edge in $E(X,V_i)$ is a bridge in both $G'$ and $G$, contradicting the assumption that $G$ is bridgeless. Hence, if $t_i=0$, then $|E(X,V_i)| \ge 3 = 3 - t_i$.
Suppose $t_i=1$ and $|E(X,V_i)|=0$ for some $i$, $1\le i\le s$. Then the unique edge in $E(G) \setminus E(G')$ that joins a vertex in $V_i$ to a vertex in $\{u,v\}$ in $G$ is a bridge in $G$, contradicting the assumption that $G$ is bridgeless. Hence, if $t_i=1$, then $|E(X,V_i)| \ge 2 = 3 - t_i$.
If $t_i = 2$, then since $|E(X,V_i)|$ and $t_i$ have different parity, we have that $|E(X,V_i)| \ge 1 = 3 - t_i$. If $t_i \ge 3$, then $|E(X,V_i)| \ge 0 \ge 3 - t_i$.~\smallqed

\medskip
\noindent
By Subclaim~2.1, $|E(X,V_i)| \ge 3-t_i$ for all $1\le i\le s$.
Recall that $\sum_{i=1}^s t_i \le 4$ and hence
\[
|E(X,\barX)| = \sum_{i=1}^s |E(X,V_i)| \ge \sum_{i=1}^s (3-t_i) = 3s - \sum_{i=1}^s t_i \ge 3\oc(G'-X)-4.
\]
This completes the proof of Claim~2.~\smallqed

\begin{unnumbered}{Claim~3}
$\oc(G'-X) \le |X|$.
\end{unnumbered}
\proof Suppose to the contrary that $\oc(G'-X) > |X|$. Since $G$ is cubic, the graph $G$, and thus the graph $G'$, has even order. Therefore, $\oc(G'-X) - |X|$ is even, implying that $\oc(G'-X) \ge |X| + 2$.  This in turn implies that
\[
|E(X,\barX)| \le 3|X| \le 3\oc(G'-X)-6,
\]
which contradicts Claim~2.~\smallqed

\medskip
By Claim~3, $\oc(G'-X) \le |X|$. By our choice of the set $X$, this implies that $\oc(G'-S) \le |S|$ for every proper subset $S \subset V(G')$. Hence by Tutte's Theorem, $G'$ has a perfect matching, $M'$ say. But then the set $M' \cup \{e\}$ is a perfect matching in $G$ containing the edge~$e$. This completes the proof of Theorem~\ref{t:matching}.~\qed

\medskip
As a consequence of Theorem~\ref{t:matching}, we have the following result.

\begin{cor}
Every $2$-connected cubic multigraph without loops is $1$-extendable.
 \label{c:matching}
\end{cor}
\proof Let $G$ be a $2$-connected cubic multigraph with no loops. We prove the result by induction on the number of parallel edges in $G$. If $G$ does not contain parallel edges, then the result follows by Theorem~\ref{t:matching}. This establishes the base case. Suppose that $G$ contains parallel edges. Then there are two vertices $u$ and $v$ in $G$ that are joined by two distinct edges, say $e$ and $e'$.

Let $D$ be a diamond with $V(D) = \{a,b,c,d\}$ and where $ab$ is the missing edge in $D$ and where $V(D) \cap V(G) = \emptyset$. Let $G'$ arises from the disjoint union, $D \cup G$, of $D$ and $G$ by removing the edge $e$ and adding the two edges $e_1 = au$ and $e_2 = bv$. Let $e_3 = cd$. By construction, the graph $G'$ is $2$-connected, cubic, without loops, and with less parallel edges than $G$. Applying the induction hypothesis to $G'$, the multigraph $G'$ is $1$-extendable.
Let $f$ be an arbitrary edge in $G$.

If the edges $e$ and $f$ are distinct,
then let $M'$ be a perfect matching in $G'$ containing $f$.
We note that $e_1\in M'$ if and only if $e_2\in M'$ if and only if $e_3\in M'$.
If $e_1\in M'$, then let $M=(M'\cap E(G))\cup\{e\}$.
Otherwise if $e_1\notin M'$, then let $M=(M'\cap E(G))$.
In both cases $M$ is a perfect matching in $G$ containing $f$.

If $e=f$,
then let $M'$ be a perfect matching in $G'$ containing $e_1$.
We note that $e_2,e_3\in M'$.
Hence $(M'\cap E(G))\cup\{e\}$ is a perfect matching in $G$ containing $e=f$.
Altogether, $G$ has a perfect matching containing $f$.~\qed



\subsection{Proof of Theorem~\ref{t:clawfree}}

We are now in a position to prove Theorem~\ref{t:clawfree}. Recall its statement.

\noindent \textbf{Theorem~\ref{t:clawfree}} \emph{If $G$ is a cubic claw-free graph of order~$n$, then $c_\ind(G) > 13n/20$.}

\medskip
\noindent \textbf{Proof of Theorem~\ref{t:clawfree}.}
Clearly, we may assume that $G$ is connected.
We proceed by induction on the order~$n$ of a connected cubic claw-free graph. If $n = 4$, then $G = K_4$ and $c_\ind(G) = 3n/4 > 13n/20$. This establishes the base case. Let $n \ge 6$ and assume that for every connected cubic claw-free graph $G'$ of order~$n'$ where $n' < n$ we have $c_\ind(G') > 13n'/20$. Let $G = (V,E)$ be a connected cubic claw-free graph of order~$n$.

By Lemma~\ref{l:triangle_diamond}, the vertex set $V$ can be uniquely partitioned into sets each of which induce a triangle or a diamond in $G$. For notational
convenience, we refer to such a partition as a
\emph{triangle}-\emph{diamond partition} of $G$, abbreviated
$\Delta$-D-partition. Every triangle and diamond induced by a set in
our $\Delta$-D-partition we call a \emph{unit} of the partition. A
unit that is a triangle we call a \emph{triangle-unit} and a unit
that is a diamond we call a \emph{diamond-unit}. We say that two
units in the $\Delta$-D-partition are \emph{adjacent} if there is an
edge joining a vertex in one unit to a vertex in the other unit.

If every unit in the $\Delta$-D-partition is a diamond-unit, then $G$
is a diamond-necklace $N_k$ with $k \ge 2$ diamonds (recall that here $n \ge 8$), and so, by Observation~\ref{necklace}, $c_\ind(G) = 3n/4 > 13n/20$ and $G \in \cN_\cub$. Hence we may assume that $G \notin \cN_\cub$. Therefore, at least one unit in the
$\Delta$-D-partition is a triangle-unit. Since every triangle-unit is
joined by three edges to vertices from other units, while every
diamond-unit is joined by two edges to vertices from other units,
there are therefore at least two triangle-units in our
$\Delta$-D-partition.

We proceed by a series of structural properties that we may assume are satisfied in the graph $G$ for otherwise the desired lower bound on $c_\ind(G)$ follows.

\begin{unnumbered}{Claim~A}
A diamond unit is adjacent to only one unit and this unit is a triangle-unit.
\end{unnumbered}
\proof Suppose first that there is a diamond unit $D$ adjacent to two distinct units $U_1$ and $U_2$. Let $V(D) = \{a,b,c,d\}$ where $ab$ is the missing edge in the diamond. Let $v_1$ and $v_2$ be the neighbors of $a$ and $b$, respectively, not in $D$. Renaming $U_1$ and $U_2$ if necessary, we may assume that $v_1 \in V(U_1)$ and $v_2 \in V(U_2)$.
We note that $v_1v_2 \notin E$. Let $G'$ be obtained from $G$ by deleting $V(D)$ and adding the edge $e = v_1v_2$. Then, $G'$ is a connected cubic claw-free graph of order~$n' = n - 4$. Applying the inductive hypothesis to $G'$, yields $c_\ind(G') > 13n'/20$.
Let $S'$ be a maximum induced $2$-regular subgraph in $G'$. 
If $e$ belongs to a cycle in $G'[S']$, then
we let $S = S' \cup \{a,b,c\}$. If $e$ does not belong to a cycle in $G'[S']$, then we may assume, renaming vertices if necessary, that $v_1 \notin S'$. In this case, we let $S = S' \cup \{a,c,d\}$. In both cases, $S$ is an induced $2$-regular subgraph in $G$ and $|S| = |S'| + 3$, implying that $c_\ind(G) \ge |S| = |S'| + 3 = c_\ind(G') + 3 > 13n'/20 + 3 = 13(n-4)/20 + 3 > 13n/20$. Hence, $D$ is adjacent to exactly one other unit. If $D$ is adjacent to another diamond unit, then $G = N_2$, a contradiction to our earlier assumption that $G \notin \cN_\cub$.~\qed

\begin{unnumbered}{Claim~B}
Two triangle-units are joined by at most two edges.
\end{unnumbered}
\proof  If two triangle-units are joined by three edges, then $G$ is a prism $C_3 \, \Box \, K_2$ of order~$n = 6$ and $c_\ind(G) = 4 > 13n/20$.~\qed

\medskip
Suppose there is a diamond-unit $D$. By Claim~A, the diamond-unit $D$ is adjacent to a triangle-unit. This triangle-unit is in turn adjacent to another triangle-unit $T$.
The subgraph of $G$ induced by these three units we call a \emph{tower}.
We call $T$  the \emph{base triangle} of the tower.
Two towers are \emph{adjacent} if there is an edge joining the two base triangles. 
A tower is shown in Figure~\ref{f:tower}.

\begin{figure}[t]
\tikzstyle{every node}=[circle, draw, fill=black!0, inner sep=0pt,minimum width=.18cm]
\begin{center}
\begin{tikzpicture}[thick,scale=.6]
  \draw(0,0) { 
    +(0.00,0.00) -- +(2.00,0.00)
    +(1.00,1.00) -- +(2.00,0.00)
    +(1.00,1.00) -- +(0.00,0.00)
    +(1.00,1.00) -- +(1.00,2.00)
    +(1.00,2.00) -- +(0.00,3.00)
    +(0.00,3.00) -- +(2.00,3.00)
    +(2.00,3.00) -- +(1.00,2.00)
    +(0.00,3.00) -- +(0.00,5.00)
    +(0.00,5.00) -- +(1.00,5.50)
    +(1.00,5.50) -- +(1.00,4.50)
    +(1.00,4.50) -- +(0.00,5.00)
    +(1.00,5.50) -- +(2.00,5.00)
    +(2.00,5.00) -- +(1.00,4.50)
    +(2.00,5.00) -- +(2.00,3.00)
    +(0.00,0.00) node{}
    +(2.00,0.00) node{}
    +(1.00,1.00) node{}
    +(1.00,2.00) node{}
    +(0.00,3.00) node{}
    +(2.00,3.00) node{}
    +(2.00,5.00) node{}
    +(0.00,5.00) node{}
    +(1.00,4.50) node{}
    +(1.00,5.50) node{}
  };
\end{tikzpicture}
\end{center}
\vskip -0.6 cm \caption{A tower.} \label{f:tower}
\end{figure}

\begin{unnumbered}{Claim~C}
If two triangle-units are joined by two edges,
then there is a third triangle-unit which is adjacent to both triangle-units.
\end{unnumbered}
\proof Suppose that there are two triangle-units $U_1$ and $U_2$ joined by two edges and there is no common adjacent (triangle) unit. For $i \in \{1,2\}$, let $u_i \in V(U_i)$ be the vertex that is not adjacent to $U_{3-i}$ and let $v_i$ be the neighbor of $u_i$ not belonging to $U_i$. Note that $v_1v_2\notin E(G)$. For $i \in \{1,2\}$, let $V(U_i) = \{u_i,w_i,x_i\}$ where $w_1w_2 \in E(G)$ and $x_1x_2 \in E(G)$. Let $G'$ be obtained from $G$ by deleting $V(U_1) \cup V(U_2)$ and adding the edge $e=v_1v_2$.
Then, $G'$ is a connected cubic claw-free graph of order~$n' = n - 6$. Applying the inductive hypothesis to $G'$, yields $c_\ind(G') > 13n'/20$.
Let $S'$ be a maximum induced $2$-regular subgraph in $G'$.
If $e$ belongs to a cycle $C$ in $G'[S']$, then let $S = S' \cup \{u_1,u_2,w_1,w_2\}$. If $e$ does not belong to a cycle in $G'[S']$, then let $S = S' \cup \{w_1,w_2,x_1,x_2\}$. In both cases, $S$ is an induced $2$-regular subgraph in $G$ and $|S| = |S'| + 4$, implying that $c_\ind(G) \ge |S| = |S'| + 4 = c_\ind(G') + 4 > 13n'/20 + 4 = 13(n-6)/20 + 4 > 13n/20$.~\qed

\begin{figure}[htb]
\tikzstyle{every node}=[circle, draw, fill=black!0, inner sep=0pt,minimum width=.182cm]
\begin{center}
\begin{tikzpicture}[thick,scale=.6]
  \draw(0,0) { 
    +(3.13,0.63) -- +(4.38,0.63)
    +(4.38,0.63) -- +(5.00,1.25)
    +(5.00,1.25) -- +(5.00,0.00)
    +(5.00,0.00) -- +(4.38,0.63)
    +(2.50,1.25) -- +(2.50,0.00)
    +(2.50,0.00) -- +(3.13,0.63)
    +(3.13,0.63) -- +(2.50,1.25)
    +(2.50,1.25) -- +(0.63,1.25)
    +(0.63,1.25) -- +(0.00,0.63)
    +(0.00,0.63) -- +(1.25,0.63)
    +(1.25,0.63) -- +(0.63,1.25)
    +(0.00,0.63) -- +(0.63,0.00)
    +(0.63,0.00) -- +(1.25,0.63)
    +(0.63,0.00) -- +(2.50,0.00)
    +(5.00,0.00) -- +(6.88,0.00)
    +(6.88,0.00) -- +(6.25,0.63)
    +(6.25,0.63) -- +(7.50,0.63)
    +(7.50,0.63) -- +(6.88,0.00)
    +(6.25,0.63) -- +(6.88,1.25)
    +(6.88,1.25) -- +(7.50,0.63)
    +(6.88,1.25) -- +(5.00,1.25)
    +(3.13,0.63) node{}
    +(4.38,0.63) node{}
    +(2.50,1.25) node{}
    +(2.50,0.00) node{}
    +(5.00,1.25) node{}
    +(5.00,0.00) node{}
    +(0.63,1.25) node{}
    +(1.25,0.63) node{}
    +(0.00,0.63) node{}
    +(0.63,0.00) node{}
    +(6.25,0.63) node{}
    +(7.50,0.63) node{}
    +(6.88,0.00) node{}
    +(6.88,1.25) node{}
  };
\end{tikzpicture}
\end{center}
\vskip -0.6 cm \caption{Two distinct non-vertex-disjoint towers with no common base triangle.} \label{f:tower2}
\end{figure}
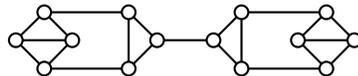

\begin{unnumbered}{Claim~D}
Every two distinct towers are vertex-disjoint.
\end{unnumbered}
\proof We show firstly that if two distinct towers are not vertex-disjoint, then they share a common base triangle. Suppose, to the contrary, that $T_1$ and $T_2$ are two distinct towers that are not vertex-disjoint but do not share a common base triangle. Then by Claim~A, the two base triangles of $T_1$ and $T_2$ cannot have exactly two vertices in common. Since $G$ is cubic, the base triangles cannot have exactly one vertex in common. Therefore the base triangles have no vertex in common and are both adjacent to a diamond-unit. The graph $G$ is determined and is shown in Figure~\ref{f:tower2}. In this case, $n = 14$ and $c_\ind(G) = 10 > 13n/20$. Hence we may assume that two distinct non-vertex-disjoint towers share a common base triangle, for otherwise the desired result follows.

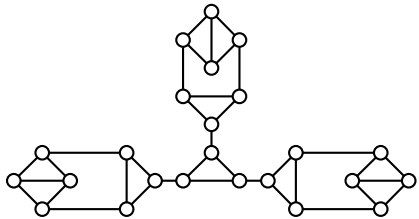
\begin{figure}[htb]
\tikzstyle{every node}=[circle, draw, fill=black!0, inner sep=0pt,minimum width=.18cm]
\begin{center}
\begin{tikzpicture}[thick,scale=.6]
  \draw(0,0) { 
    +(5.63,0.63) -- +(6.25,1.25)
    +(6.25,1.25) -- +(6.25,0.00)
    +(6.25,0.00) -- +(5.63,0.63)
    +(2.50,1.25) -- +(2.50,0.00)
    +(2.50,0.00) -- +(3.13,0.63)
    +(3.13,0.63) -- +(2.50,1.25)
    +(2.50,1.25) -- +(0.63,1.25)
    +(0.63,1.25) -- +(0.00,0.63)
    +(0.00,0.63) -- +(1.25,0.63)
    +(1.25,0.63) -- +(0.63,1.25)
    +(0.00,0.63) -- +(0.63,0.00)
    +(0.63,0.00) -- +(1.25,0.63)
    +(0.63,0.00) -- +(2.50,0.00)
    +(6.25,0.00) -- +(8.13,0.00)
    +(8.13,0.00) -- +(7.50,0.63)
    +(7.50,0.63) -- +(8.75,0.63)
    +(8.75,0.63) -- +(8.13,0.00)
    +(7.50,0.63) -- +(8.13,1.25)
    +(8.13,1.25) -- +(8.75,0.63)
    +(8.13,1.25) -- +(6.25,1.25)
    +(4.38,1.25) -- +(3.75,0.63)
    +(3.75,0.63) -- +(5.00,0.63)
    +(5.00,0.63) -- +(4.38,1.25)
    +(3.13,0.63) -- +(3.75,0.63)
    +(5.00,0.63) -- +(5.63,0.63)
    +(4.38,4.38) -- +(4.38,3.13)
    +(3.75,3.75) -- +(3.75,2.50)
    +(5.00,3.75) -- +(5.00,2.50)
    +(5.00,3.75) -- +(4.38,3.13)
    +(4.38,3.13) -- +(3.75,3.75)
    +(3.75,3.75) -- +(4.38,4.38)
    +(4.38,4.38) -- +(5.00,3.75)
    +(5.00,2.50) -- +(3.75,2.50)
    +(3.75,2.50) -- +(4.38,1.88)
    +(4.38,1.88) -- +(5.00,2.50)
    +(4.38,1.88) -- +(4.38,1.25)
    +(3.13,0.63) node{}
    +(5.63,0.63) node{}
    +(2.50,1.25) node{}
    +(2.50,0.00) node{}
    +(6.25,1.25) node{}
    +(6.25,0.00) node{}
    +(0.63,1.25) node{}
    +(1.25,0.63) node{}
    +(0.00,0.63) node{}
    +(0.63,0.00) node{}
    +(7.50,0.63) node{}
    +(8.75,0.63) node{}
    +(8.13,0.00) node{}
    +(8.13,1.25) node{}
    +(3.75,0.63) node{}
    +(4.38,1.25) node{}
    +(5.00,0.63) node{}
    +(4.38,1.88) node{}
    +(3.75,2.50) node{}
    +(5.00,2.50) node{}
    +(3.75,3.75) node{}
    +(4.38,4.38) node{}
    +(5.00,3.75) node{}
    +(4.38,3.13) node{}
  };
\end{tikzpicture}
\end{center}
\vskip -0.6 cm \caption{Three towers sharing a common base triangle.} \label{f:tower3}
\end{figure}

%
Suppose that three towers share a common base triangle. Then the graph $G$ is determined and is shown in Figure~\ref{f:tower3}. In this case, $n = 24$ and $c_\ind(G) \ge 18 > 13n/20$. Hence we may assume that no three towers share a common base triangle, for otherwise the desired result follows.

Suppose $G_1$ and $G_2$ are two distinct towers that share a common base triangle $T$. Let $v$ be the vertex of $T$ that has a neighbor in $G-V(G_1)-V(G_2)$. Let $r_1$ be the neighbor of $v$ not in $G_1$. By Claim~A, $r_1$ belongs to a triangle-unit, say $T^*$. Let $V(T^*) = \{r_1,r_2,r_3\}$. Let $s_2$ and $s_3$ be the neighbors of $r_2$ and $r_3$, respectively, not in $T^*$. By our earlier assumption, no three towers share a common triangle. Hence by Claim~C, $s_2$ and $s_3$ belong to different units. In particular, we note that $s_2s_3\notin E(G)$.

Let $G'$ be the graph obtained from $G$ by deleting the vertices in $V(G_1) \cup V(G_2) \cup V(T^*)$ and adding the edge $s_2s_3$. Then, $G'$ is a connected cubic claw-free graph of order $n' = n - 20$. Let $C_1$ and $C_2$ be induced $5$-cycles in $G_1$ and $G_2$, respectively. Applying the inductive hypothesis to $G'$, we have $c_\ind(G') > 13n'/20$.
Let $S'$ be a maximum induced $2$-regular subgraph in $G'$. 
If $\{s_2,s_3\} \subset S'$, let $S = S' \cup \{r_2,r_3\} \cup V(C_1) \cup V(C_2) \cup V(T)$. If $|S' \cap \{s_2,s_3\}| \le 1$, let $S = S' \cup V(C_1) \cup V(C_2) \cup V(T)$. In both cases, $S$ is an induced $2$-regular subgraph in $G$, and so $c_\ind(G) \ge |S| \ge |S'| + 13 = c_\ind(G') + 13 > 13n'/20 + 13 = 13(n-20)/20 + 13 = 13n/20$.~\qed

\medskip
We now return to the proof of Theorem~\ref{t:clawfree}. By our assumptions to date, the connected claw-free graph $G$ satisfies Claims A-D.
We now construct an induced $2$-regular subgraph in $G$ explicitly. In order to do this, we introduce several auxiliary graphs. By Claim~D, every two distinct towers are vertex-disjoint. Let $H$ arise from $G$ as follows: First we replace every tower by a degree-$2$ vertex. Thereafter we replace every triangle-unit in the resulting reduced graph, by a degree-$3$ vertex.
Let $h$ be the bijective function that maps every degree-$2$ vertex of $H$ to its corresponding tower in $G$ and
every degree-$3$ vertex of $H$ to its corresponding triangle in $G$.
Note that $H$ has no loops but parallel edges are possible.

We color in $H$ every degree-$2$ vertex red and every degree-$3$ vertex which is incident to a bridge yellow. Let $H_1, \ldots, H_k$ be the non-trivial and non-cycle blocks of $H$. Then every uncolored vertex belongs to $H_i$ for some $i$, $1 \le i \le k$, has degree~$3$, and all its three incident edges belong to the same non-trivial block of $H$. Note that in $H_1,\ldots,H_k$ exactly the colored vertices have degree~$2$. If $u$ and $v$ are two distinct uncolored vertices in $H_i$ for some $i$, $1 \le i \le k$, that are joined by a path $P$ every internal vertex of which is a colored vertex (colored red or yellow), then we refer to $P$ as a \emph{colored} $u$-$v$-path.

For every $1\leq i \leq k$, let $B_i$ be the graph that arises from the subgraph of $H$ induced by all uncolored vertices of $H_i$ by adding, for each colored $u$-$v$-path between two uncolored vertices $u$ and $v$ in $H_i$, an edge between $u$ and $v$. We say that such an added edge is colored.
Note that every $1\leq i \leq k$, the graphs $B_i$ may have parallel edges.
However, $B_i$ does not contain loops because $H_i$ is neither an edge nor a cycle.
Moreover, we also use the function $h$ to map every vertex of $B_i$ to its corresponding triangle in $G$.

In the following, we take a set of cycles in $H$ or $B_1,\ldots, B_k$ and construct from this set one part of a $2$-regular subgraph of $G$. In order to do this, we introduce the notion of \textit{lifting} a cycle.
Let $C$ be a cycle in $B_i$. Let $C'$ be the cycle of $H_i$ that arises from $C$ by replacing each colored edge with its corresponding colored path. Let $C''$ be the cycle of $G$ that uses exactly two vertices of all (triangle-) units in $h(V(C'))$ and respects the naturally given order of $C'$.
We say that $C$ (or $C'$) is \textit{lifted} to $C''$ and write $\ell(C)= C''$.
By abusing the terminology, we also write $\ell(C')=C''$. We remark that $C''$ is an induced cycle of $G$ because if $e \in E(G) \setminus E(C'')$ and $e = xy$ is incident with a vertex in $C''$, then $e$ belongs to a triangle-unit $T_e$, say, where $V(T_e) = \{x,y,z\}$ but there are only two vertices in $T_e$ that belong to the cycle $C''$, namely $y$ and $z$.
Moreover, we extend the idea of lifting a cycle to collections of (pairwise disjoint) cycles by lifting each cycle separately. Note that if we lift a collection of pairwise disjoint cycles $\cD$ , then $\ell(\cD)$ is a $2$-regular induced subgraph of $G$.

One core idea of our proof is to remove from each $B_i$, $1 \le i \le k$, a perfect matching $M_i$ resulting in a collection $B_i-M_i$  of cycles which we then lift to induced cycles in $G$.
Applying Corollary~\ref{c:matching} to the $2$-connected cubic multigraph $B_i$,
the multigraph $B_i$ is $1$-extendable and therefore we have the flexibility to force any specific edge $e$ of $B_i$ to be in the perfect matching $M_i$ or not.
In order to force $e\notin M_i$ we force a adjacent edge of $e$ to be contained in $M_i$.

Before we are able to construct a large $2$-regular induced subgraph of $G$, we describe some properties of the components $H^c_1,\ldots,H^c_{\ell}$ induced by the colored vertices of $H$.
It is easy to see that $H^c_{i}$ is a cactus; that is, every block in $H^c_{i}$ is either a trivial block or a cycle block. By the definition of a block, the components $H^c_1,\ldots,H^c_{\ell}$ and the non-trivial and non-cycle blocks $H_1,\ldots H_k$ can be structured in a tree $T$ with vertex set $\{t_{H_1},\ldots t_{H_k},t_{H^c_{1}},\ldots,t_{H^c_{\ell}}\}$, where $t_{H_i}$ and $t_{H^c_{j}}$ are joined by an edge if $V(H_i) \cap V(H^c_{j}) \ne  \emptyset$. Note that for all distinct $i,j$, we have $V(H_i) \cap V(H_j) = \emptyset$ and $V(H^c_{i}) \cap V(H^c_{j}) = \emptyset$.
We call $B_{1},\ldots,B_{k}$ and $H^c_{1},\ldots,H^c_{\ell}$ the pieces of $G$.
This implies that $V(B_{1}),\ldots,V(B_{k}),V(H^c_{1}),\ldots,V(H^c_{\ell})$ is a partition of $V(H)$.

In order to construct a collection of induced cycles $\cC$ of $G$, in other words a $2$-regular subgraph of $G$, our strategy is as follows: Initially, we let $\cC=\emptyset$. Suppose $k\ge 1$;  we treat the case $k=0$ in the very end of the proof. Let $T$ be  rooted  at $t_{H_1}$. For $i \ge 0$, we define the $i$th level of $T$ to be the set of all vertices in $T$ at distance~$i$ from the root, $t_{H_1}$, of $T$. In particular, level~$0$ consists of the root $t_{H_1}$ of $T$. We start with some perfect matching $M_1$ of $B_1$ and add $\ell(B_1-M_1)$ to $\cC$.
Thus $|\cC \cap h(B_1)| \ge \frac{2}{2}n(h(B_1))$.
We then consider the levels of $T$ in turn, starting with the vertices in the 1st level, proceeding to vertices in the 2nd level, and so on. Hence after $\ell(B_1-M_1)$ is added to $\cC$, as our next step we consider all pieces $H^c_{i}$ such that $t_{H^c_{i}}$ is adjacent to $t_{H_1}$. Thereafter, we consider all pieces $H_j$ such that $t_{H_j}$ is at distance~$2$ from $t_{H_1}$ in $T$ and we continue this process until all vertices of $T$ have been examined. Since $T$ is a tree and the influence of our construction points away from $t_{H_1}$, it is sufficient to explain our construction for some piece $H_i$ and for some piece $H^c_{j}$ if its parent vertex in $T$ has already been examined, respectively; in other words we proceed by induction on the levels of the vertices in $T$.

Let $d\ge 1$ and for the purpose of induction,
we assume that $|\cC \cap h(P)| \ge \frac{13}{20}n(h(P))$ for all pieces $P$ such that $t_P$ is at level less than~$d$ in $T$, with strict inequality for $P=B_1$
where $\cC$ is a $2$-regular induced subgraph of $G$.
Let $t_P$ be a vertex at level~$d$ in $T$.

Suppose that $P = H^c_{j}$ for some $1\leq j \leq \ell$. Let $t_{H_i}$ be the parent vertex of $t_{H^c_{j}}$. For all towers in $h(H^c_{j})$, add one of their two induced $5$-cycles to $\cC$.

\begin{unnumbered}{Claim~E}
If $H^c_{j}$ contains no yellow vertex, then $\cC$ can be extended to a $2$-regular graph of $G$ by adding to it at least~$\frac{13}{20}n(h(H^c_{j}))$ vertices.
\end{unnumbered}
\proof Suppose that $H^c_{j}$ is a path on $p$ red colored vertices. In this case, we note that $t_{H^c_{j}}$ is a leaf in $T$. If $M_i$ contains the colored edge in $B_i$ that is associated with $H^c_{j}$, then there is no cycle in $\cC$ that contains a vertex of $h(H^c_{j})$.
Further no vertex from a cycle in $\cC$ is adjacent to a vertex in $h(H^c_{j})$.
We choose the first base triangle, and thereafter every second base triangle, of $H^c_{j}$ to $\cC$. In this way, $\lceil p/2 \rceil$ base triangles of $H^c_{j}$ are added to $\cC$, implying that $\cC$ is a $2$-regular graph of $G$. Further since $n(h(H^c_{j})) = 10p$, we have
\[
|\cC \cap h(H^c_{j})| =
 5p + 3 \left\lceil \frac{p}{2} \right\rceil  \ge
 \frac{13p}{2} = \frac{13}{20}n(h(H^c_{j})).
\]
If $M_i$ does not contain the colored edge that is associated to $H^c_{j}$, then there is a cycle in $\cC$ that contains exactly two vertices of every base triangle from each tower of $h(H^c_{j})$. In this case, $|\cC \cap h(H^c_{j})| = 7p > \frac{13}{20}n(h(H^c_{j}))$.~\smallqed

\medskip
In view of Claim~E, we may assume that $H^c_{j}$ contains at least one yellow vertex. For every cycle $C$ in $H^c_{j}$, add the lifted cycle $\ell(C)$ to $\cC$. Note that $\cC$ is still a $2$-regular subgraph of $G$ and $\frac{2}{3}n(\ell(C))$ vertices have been added to $\cC$. Trivially, all vertices of $H^c_{j}$ that do not belong to a cycle of $H^c_{j}$ induce a forest $F$. Since $F$ is bipartite, we can partition $F$ into two independent sets $I_1$ and $I_2$. Both of them represent a collection of induced cycles (in fact, triangles) in $G$.
We choose later one of these collections and add it to $\cC$ and show that at least one choice is good enough for our purpose by showing that the average of our two choices is already good enough.

Let $v\in V(H^c_{j})\cap I_p$ for some $p\in \{1,2\}$ and suppose we choose $I_p$. Then we do the following: if $v$ is a yellow colored vertex, then add $h(v)$ to $\cC$ and if $v$ is a red colored vertex, then add the base triangle of $h(v)$ to $\cC$.

Note that, if $v\in V(H^c_{j})$ is a vertex that belongs to $H_{i'}$ for some $i'\ne i$, then $t_{H_{i'}}$ is a child vertex of $t_{H^c_{j}}$ and all vertices of $H^c_{j}$ which belong to $H_{i'}$ belong to the same colored edge $e_{i'}$ of $B_{i'}$. By Corollary~\ref{c:matching}, we are able to decide whether $e_{i'}$ is in the perfect matching $M_{i'}$ of $B_{i'}$ or not.
That means, if we decide that $e_{i'}\in M_{i'}$, then after later adding $\ell(B_{i'}-M_{i'})$ to $\cC$, there is no cycle in $\cC$ completely lying in $H_{i'}$ and using vertices of $G$ which correspond to the colored edge $e_{i'}$.
In the latter case, there is such a cycle using two vertices of each triangle and each base triangle of a tower which correspond to the colored edge $e_{i'}$. Thus we have more freedom than only choosing one of the independent sets $I_1$ or $I_2$.

Next, we analyze the average contribution of $|\cC \cap h(v)|$ to $|\cC \cap h(H^c_{j})|$ for $v\in V(H^c_{j})$. For simplification, we consider the average contribution
\[
a(v)=|\cC \cap h(v)|- \frac{13}{20}n(h(v))
\]
and show that
\[
\sum_{v\in V(H^c_{j})}a(v) \ge 0.
\]

Let $v \in V(H^c_{j})$. We distinguish several cases:
\begin{enumerate}

\item The vertex $v$ belonging to $H_i$ is a yellow-colored vertex. If there is already a cycle $C\in\cC$ with $h(v) \cap C$ that belongs completely to $H_i$, then this cycle uses two vertices of $h(v)$ and we remove $v$ from $I_1 \cup I_2$. This implies that $a(v)= 2- \frac{39}{20}= \frac{1}{20}$. If this is not the case, then $v$ is either in $I_1$ or $I_2$ and hence $a(v)= \frac{1}{2}(0+3)- \frac{39}{20}= -\frac{9}{20}$. Thus, $a(v)\ge -\frac{9}{20}$.

\item The vertex $v$ belonging to $H_i$ is a red-colored vertex. If there is already a cycle in $C\in\cC$ with $h(v)\cap C\ne \emptyset$ that belongs completely to $H_i$, then this cycle uses two vertices of $h(v)$ and we remove $v$ from $I_1 \cup I_2$. This implies that $a(v) = 7- \frac{130}{20}= \frac{1}{2}$. If this is not the case, then $v$ is either in $I_1$ or $I_2$ and hence $a(v)= \frac{1}{2}(5+8)- \frac{130}{20}= 0$. Thus, $a(v)\ge 0$.

\item The vertex $v$ is yellow-colored and does not belong to any $H_{i'}$ and no colored cycle. Since $v$ is in exactly one independent set, we obtain $a(v)= \frac{1}{2}(0+3) - \frac{39}{20}= - \frac{9}{20}$.
	
\item The vertex $v$ is red-colored and does not belong to any $H_{i'}$ and no colored cycle. Since $v$ is in exactly one independent set, we obtain $a(v)= \frac{1}{2}(5+8) - \frac{130}{20}= 0$.

\item The vertex $v$ belongs to a colored cycle, $C$, of length~$r$. We calculate $\sum_{v\in V(C)}a(v)$.
	
\begin{enumerate}
	  \item Suppose first that $C$ contains only one yellow vertex. By Claim~C, the cycle $C$ contains at least two red vertices. Hence \[
\sum_{v\in V(C)} a(v) = (7(r-1)+2)-(10(r-1)+3) \cdot \frac{13}{20}\ge  \frac{21}{20}.
\]
				
\item Suppose now that $C$ contains at least two yellow vertices. Since $|\cC \cap h(v)|\ge \frac{2}{3} n(h(v))$ for every $v\in V(C)$, we conclude $\sum_{v\in V(C)}a(v)\ge 0$.
\end{enumerate}

\item The vertex $v$ belongs to some $H_{i'}$ where $i' \ne i$. Let $e_{i'}$ be a colored edge belonging to $B_{i'}$ that is associated with $H^c_{j}$. Let $P$ be the colored path in $H$ that corresponds to the colored edge $e_{i'}$.
	
\begin{enumerate}
\item If $P$ contains at least two yellow vertices, then remove all vertices of $P$ from $I_1$ and $I_2$ and we decide that $e_{i'}\notin M_{i'}$.   This implies that $|\cC \cap h(v)|\ge \frac{2}{3} n(h(v))$ for every $v\in V(C)$, and therefore that $\sum_{v\in V(C)}a(v)\ge 0$.

\item If $P$ contains contains exactly one yellow vertex and at least one red vertex, then remove all vertices of $P$ from $I_1$ and $I_2$ and we decide that $e_{i'}\notin M_{i'}$. If $v$ is the yellow vertex of $P$, then $a(v) = \frac{1}{20}$ while for every red vertex $w$ of $P$, we have $a(w) = 7 - \frac{13}{2} = \frac{1}{2}$. Thus since $P$ has at least one red vertex,
    \[ \sum_{v\in V(P)}a(v) \ge \frac{1}{20}+ \frac{1}{2}= \frac{11}{20}.\]
				
\item Suppose $P$ contains exactly one vertex $v$. The vertex $v$ is necessarily a yellow vertex. Let $p$ be such that $v \in I_p$. If we later choose $I_p$, then we decide that $e_{i'}\in M_{i'}$. Otherwise, if $I_p$ is not chosen, we decide that $e_{i'}\notin M_{i'}$ and we remove $v$ from $I_p$. In the first case we have $|\cC \cap h(v)|=3$ and in the latter case $|\cC \cap h(v)|=2$. 		Thus, $a(v)= \frac{1}{2}(3+2)- \frac{39}{20}= \frac{11}{20}$.
	\end{enumerate}
\end{enumerate}

We remark that if $v \in V(H^c_{j})$, then $a(v)$ is only negative in Case~1 and Case~3.

\begin{unnumbered}{Claim~F}
The cases 1 and 3 occur at most as often as the cases 5(a), 6(b) and 6(c).
\end{unnumbered}
\proof
In the cactus $H^c_{j}$ we contract every cycle to a vertex as well as all vertices which belong to some $H_{i'}$ for $i'\ne  i$, respectively.
Denote this graph by $T^j$. Note that $T^j$ is a tree.
Let $e_i$ be the colored edge in $B_i$ that is associated with $H^c_{j}$ and let $p$ be the number of yellow vertices in this path.
Partition $T^j$ in $p$ subtrees such that each tree contains exactly one of these $p$ yellow vertices. Contract all vertices that correspond to $e_i$ in each subtree to one vertex and denote the resulting trees by $T_1, \ldots, T_p$. 
Let $1\leq q \leq p$.
Each leaf of $T_q$, beside the leaf which corresponds to case~1, corresponds to one occurrence of either case 5(a), 6(b) and 6(c).
Furthermore, every occurrence of case 3 guarantees one vertex of degree~$3$ in $T_q$.
Since every tree has at least two more leaves than the number of degree~$3$ vertices,
we conclude the desired result.~\qed

\medskip
Since the positive contribution of each of the cases 5(a), 6(b) and 6(c) is as least as large as the absolute value of the negative contribution of each of the cases 1 and 3, we conclude that \[
\sum_{v\in V(H^c_{j})}a(v) \ge 0.
\]

This implies that $I_1$ or $I_2$ is a good choice for us and we accordingly add all associated triangles to our choice to $\cC$. This implies that
\[
|\cC \cap h(H^c_{j})| \ge \frac{13}{20} n(h(H^c_{j})).
\]

Suppose next that $P = B_i$ for some $2\leq i \in k$. Let $t_{H^c_{j}}$ be the parent vertex of $t_{H_i}$ for some $1\leq j \leq \ell$.
Let $M_i$ be a perfect matching of $B_i$.
By our earlier discussion, there is an colored edge $e_{i'}$ for which we decided whether $e_{i'} \in M_{i'}$ or not.
Recall that $B_{i'}$ is a $2$-connected multigraph.
By Corollary \ref{c:matching},
we obtain the existence of a perfect matching $M_{i'}$ of $B_{i'}$ which contains $e_{i'}$ or not accordingly to our previous decision.
With this choice of $M_i$, we extend $\cC$ to a larger $2$-regular induced subgraph of $G$ by adding  $\ell(B_i-M_i)$ to $\cC$. This implies that
\[
|\cC \cap h(B_i)| = \frac{2}{3}n(h(B_i) > \frac{13}{20}n(h(B_i)).
\]
Since $V(B_{1}),\ldots,V(B_{k}),V(H^c_{1}),\ldots,V(H^c_{\ell})$
is a partition of $V(H)$, it remains for us to consider the case $k=0$. However, this case follows by the same approach as before except that here only the cases 3, 4 and 5 can occur. This completes the proof of Theorem~\ref{t:clawfree}.~\qed

\medskip
That the bound in Theorem~\ref{t:clawfree} is asymptotically tight may be seen as follows.

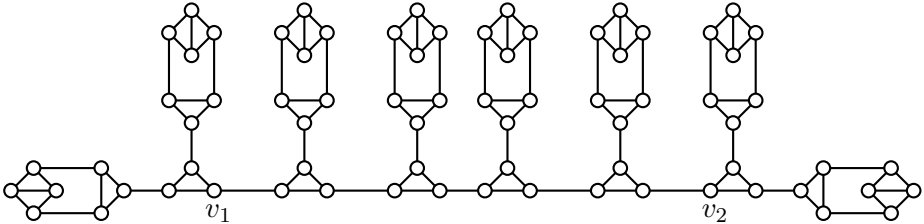
\begin{figure}[htb]
\tikzstyle{every node}=[circle, draw, fill=black!0, inner sep=0pt,minimum width=.18cm]
\begin{center}
\begin{tikzpicture}[thick,scale=1.2]
  \draw(0,0) { 
    +(6.75,1.00) -- +(6.75,0.50)
    +(6.50,1.25) -- +(7.00,1.25)
    +(7.00,1.25) -- +(6.75,1.00)
    +(6.75,1.00) -- +(6.50,1.25)
    +(6.50,1.25) -- +(6.50,2.00)
    +(6.50,2.00) -- +(6.75,2.25)
    +(6.75,2.25) -- +(7.00,2.00)
    +(7.00,2.00) -- +(6.75,1.75)
    +(6.75,1.75) -- +(6.50,2.00)
    +(6.75,2.25) -- +(6.75,1.75)
    +(7.00,2.00) -- +(7.00,1.25)
    +(6.75,0.50) -- +(6.50,0.25)
    +(6.50,0.25) -- +(7.00,0.25)
    +(7.00,0.25) -- +(6.75,0.50)
    +(7.00,0.25) -- +(7.75,0.25)
    +(7.75,0.25) -- +(8.00,0.50)
    +(8.25,0.25) -- +(8.00,0.50)
    +(8.25,0.25) -- +(7.75,0.25)
    +(5.75,0.25) -- +(6.50,0.25)
    +(5.75,0.25) -- +(5.25,0.25)
    +(5.25,0.25) -- +(4.75,0.25)
    +(4.75,0.25) -- +(4.25,0.25)
    +(4.25,0.25) -- +(3.50,0.25)
    +(3.50,0.25) -- +(3.00,0.25)
    +(3.00,0.25) -- +(2.25,0.25)
    +(2.25,0.25) -- +(2.00,0.50)
    +(2.00,0.50) -- +(1.75,0.25)
    +(1.75,0.25) -- +(2.25,0.25)
    +(1.75,0.25) -- +(1.25,0.25)
    +(2.00,0.50) -- +(2.00,1.00)
    +(2.00,1.00) -- +(1.75,1.25)
    +(1.75,1.25) -- +(2.25,1.25)
    +(2.25,1.25) -- +(2.00,1.00)
    +(1.75,1.25) -- +(1.75,2.00)
    +(1.75,2.00) -- +(2.00,2.25)
    +(2.00,2.25) -- +(2.00,1.75)
    +(2.00,1.75) -- +(1.75,2.00)
    +(2.00,2.25) -- +(2.25,2.00)
    +(2.25,2.00) -- +(2.00,1.75)
    +(2.25,2.00) -- +(2.25,1.25)
    +(1.00,0.50) -- +(1.25,0.25)
    +(1.25,0.25) -- +(1.00,0.00)
    +(1.00,0.00) -- +(1.00,0.50)
    +(1.00,0.50) -- +(0.25,0.50)
    +(0.25,0.50) -- +(0.00,0.25)
    +(0.00,0.25) -- +(0.25,0.00)
    +(0.25,0.00) -- +(0.50,0.25)
    +(0.50,0.25) -- +(0.25,0.50)
    +(0.00,0.25) -- +(0.50,0.25)
    +(0.25,0.00) -- +(1.00,0.00)
    +(9.75,0.50) -- +(10.00,0.25)
    +(10.00,0.25) -- +(9.75,0.00)
    +(9.75,0.00) -- +(9.50,0.25)
    +(9.50,0.25) -- +(9.75,0.50)
    +(9.50,0.25) -- +(10.00,0.25)
    +(9.75,0.50) -- +(9.00,0.50)
    +(9.00,0.50) -- +(9.00,0.00)
    +(9.00,0.00) -- +(9.75,0.00)
    +(9.00,0.50) -- +(8.75,0.25)
    +(8.75,0.25) -- +(9.00,0.00)
    +(8.75,0.25) -- +(8.25,0.25)
    +(8.00,0.50) -- +(8.00,1.00)
    +(8.00,1.00) -- +(7.75,1.25)
    +(7.75,1.25) -- +(8.25,1.25)
    +(8.25,1.25) -- +(8.00,1.00)
    +(7.75,1.25) -- +(7.75,2.00)
    +(7.75,2.00) -- +(8.00,2.25)
    +(8.00,2.25) -- +(8.00,1.75)
    +(8.00,1.75) -- +(7.75,2.00)
    +(8.00,2.25) -- +(8.25,2.00)
    +(8.25,2.00) -- +(8.00,1.75)
    +(8.25,2.00) -- +(8.25,1.25)
    +(5.50,1.00) -- +(5.50,0.50)
    +(5.50,0.50) -- +(5.75,0.25)
    +(5.50,0.50) -- +(5.25,0.25)
    +(5.50,1.00) -- +(5.25,1.25)
    +(5.25,1.25) -- +(5.75,1.25)
    +(5.75,1.25) -- +(5.50,1.00)
    +(5.25,1.25) -- +(5.25,2.00)
    +(5.25,2.00) -- +(5.50,2.25)
    +(5.50,2.25) -- +(5.50,1.75)
    +(5.50,1.75) -- +(5.25,2.00)
    +(5.50,2.25) -- +(5.75,2.00)
    +(5.75,2.00) -- +(5.50,1.75)
    +(5.75,2.00) -- +(5.75,1.25)
    +(4.50,0.50) -- +(4.25,0.25)
    +(4.75,0.25) -- +(4.50,0.50)
    +(4.50,0.50) -- +(4.50,1.00)
    +(4.50,1.00) -- +(4.25,1.25)
    +(4.25,1.25) -- +(4.75,1.25)
    +(4.75,1.25) -- +(4.50,1.00)
    +(4.25,1.25) -- +(4.25,2.00)
    +(4.25,2.00) -- +(4.50,2.25)
    +(4.50,2.25) -- +(4.50,1.75)
    +(4.50,1.75) -- +(4.25,2.00)
    +(4.50,2.25) -- +(4.75,2.00)
    +(4.75,2.00) -- +(4.50,1.75)
    +(4.75,2.00) -- +(4.75,1.25)
    +(3.50,0.25) -- +(3.25,0.50)
    +(3.25,0.50) -- +(3.00,0.25)
    +(3.25,0.50) -- +(3.25,1.00)
    +(3.25,1.00) -- +(3.00,1.25)
    +(3.00,1.25) -- +(3.50,1.25)
    +(3.50,1.25) -- +(3.50,2.00)
    +(3.50,1.25) -- +(3.25,1.00)
    +(3.00,1.25) -- +(3.00,2.00)
    +(3.00,2.00) -- +(3.25,2.25)
    +(3.25,2.25) -- +(3.50,2.00)
    +(3.50,2.00) -- +(3.25,1.75)
    +(3.25,1.75) -- +(3.00,2.00)
    +(3.25,2.25) -- +(3.25,1.75)
    +(0.00,0.25) node{}
    +(0.25,0.50) node{}
    +(0.25,0.00) node{}
    +(0.50,0.25) node{}
    +(1.00,0.50) node{}
    +(1.00,0.00) node{}
    +(1.25,0.25) node{}
    +(1.75,0.25) node{}
    +(2.00,0.50) node{}
    +(2.25,0.25) node{}
    +(2.00,1.00) node{}
    +(1.75,1.25) node{}
    +(2.25,1.25) node{}
    +(1.75,2.00) node{}
    +(2.25,2.00) node{}
    +(2.00,1.75) node{}
    +(3.00,0.25) node{} +(2.30,0) node[rectangle, draw=white!0, fill=white!100]{$v_1$}
    +(3.50,0.25) node{}
    +(3.25,0.50) node{}
    +(3.25,1.00) node{}
    +(3.00,1.25) node{}
    +(3.50,1.25) node{}
    +(2.00,2.25) node{}
    +(3.25,1.75) node{}
    +(3.25,2.25) node{}
    +(3.00,2.00) node{}
    +(3.50,2.00) node{}
    +(4.25,0.25) node{}
    +(4.75,0.25) node{}
    +(4.50,0.50) node{}
    +(4.50,1.00) node{}
    +(4.25,1.25) node{}
    +(4.75,1.25) node{}
    +(4.25,2.00) node{}
    +(4.50,2.25) node{}
    +(4.75,2.00) node{}
    +(4.50,1.75) node{}
    +(5.25,0.25) node{}
    +(5.75,0.25) node{}
    +(5.50,0.50) node{}
    +(5.50,1.00) node{}
    +(5.25,1.25) node{}
    +(5.75,1.25) node{}
    +(5.50,1.75) node{}
    +(5.25,2.00) node{}
    +(5.50,2.25) node{}
    +(5.75,2.00) node{}
    +(6.50,0.25) node{}
    +(7.00,0.25) node{}
    +(7.75,0.25) node{}
    +(8.00,0.50) node{}
    +(8.25,0.25) node{} +(7.80,0) node[rectangle, draw=white!0, fill=white!100]{$v_2$}
    +(6.75,0.50) node{}
    +(8.75,0.25) node{}
    +(9.00,0.50) node{}
    +(9.00,0.00) node{}
    +(9.75,0.50) node{}
    +(10.00,0.25) node{}
    +(9.50,0.25) node{}
    +(9.75,0.00) node{}
    +(6.75,1.00) node{}
    +(6.50,1.25) node{}
    +(7.00,1.25) node{}
    +(7.00,2.00) node{}
    +(6.50,2.00) node{}
    +(6.75,2.25) node{}
    +(6.75,1.75) node{}
    +(8.00,1.00) node{}
    +(7.75,1.25) node{}
    +(8.25,1.25) node{}
    +(8.00,2.25) node{}
    +(7.75,2.00) node{}
    +(8.25,2.00) node{}
    +(8.00,1.75) node{}
  };
\end{tikzpicture}
\end{center}
\vskip -0.6 cm \caption{The cubic claw-free graph $G_2$.} \label{f:G2}
\end{figure}

\begin{prop}
For any given $\epsilon > 0$, there exists a connected cubic claw-free graph of sufficiently large order $n(G)$ such that
\[
\frac{c_\ind(G)}{n(G)} < \frac{13}{20} + \epsilon.
\]
 \label{p:tight}
\end{prop}
\proof For a given (fixed) $\epsilon > 0$, let $k$ be a positive integer satisfying
\[
k > \frac{18}{400 \epsilon} - \frac{17}{10}.
\]
Let $G_k$ be the connected cubic claw-free graph obtained as follows. Let $T$ be two towers that share a common base triangle, and let $T_1$ and $T_2$ be two vertex disjoint copies of $T$.
Let $v_1$ and $v_2$ be the vertices of degree~$2$ in $T_1$ and $T_2$, respectively.
Let $G_k$ be obtained from $T_1 \cup T_2$ by joining $v_1$ and $v_2$ with a path with $2k$ internal vertices and then replacing each of the $2k$ internal vertices  (of degree~$2$) on this path with a tower in such a way that the resulting graph is cubic.
The graph $G_2$, for example, is illustrated in Figure~\ref{f:G2}. Then, $G_k$ has order~$n(G) = 20k + 34$ and $c_\ind(G) = 13k + 23$. Thus, by our choice of $k$,
\[
\frac{c_\ind(G)}{n(G)} = \frac{13k + 23}{20k + 34} < \frac{13}{20} + \epsilon.
\]
This completes the proof of Proposition~\ref{p:tight}.~\qed


\section{Closing Conjectures}

We believe that the lower bound established in  Corollary~\ref{c:cubic_bd} on the induced cycle number of a cubic graph is not optimal and
conjecture the following stronger result.

\begin{conj}
If $G$ is a cubic graph of order~$n$, then $c_\ind(G) \ge n/2$.
 \label{conj_cubic}
\end{conj}

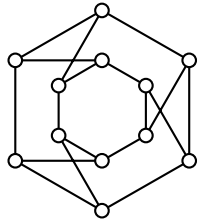
\begin{figure}[htb]
\tikzstyle{every node}=[circle, draw, fill=black!0, inner sep=0pt,minimum width=.18cm]
\begin{center}
\begin{tikzpicture}[thick,scale=.8]
  \draw(0,0) { 
    +(1.44,2.50) -- +(0.72,2.08)
    +(0.72,1.25) -- +(0.72,2.08)
    +(1.44,0.83) -- +(0.72,1.25)
    +(1.44,0.83) -- +(2.17,1.25)
    +(2.17,1.25) -- +(2.17,2.08)
    +(2.17,2.08) -- +(1.44,2.50)
    +(1.44,3.33) -- +(0.00,2.50)
    +(0.00,2.50) -- +(0.00,0.83)
    +(0.00,0.83) -- +(1.44,0.00)
    +(1.44,0.00) -- +(2.89,0.83)
    +(2.89,0.83) -- +(2.89,2.50)
    +(2.89,2.50) -- +(1.44,3.33)
    +(1.44,3.33) -- +(0.72,2.08)
    +(0.00,2.50) -- +(1.44,2.50)
    +(0.72,1.25) -- +(1.44,0.00)
    +(1.44,0.83) -- +(0.00,0.83)
    +(2.17,2.08) -- +(2.89,0.83)
    +(2.17,1.25) -- +(2.89,2.50)
    +(1.44,2.50) node{} 
    +(0.72,2.08) node{}
    +(0.72,1.25) node{}
    +(1.44,0.83) node{}
    +(2.17,1.25) node{}
    +(2.17,2.08) node{}
    +(1.44,3.33) node{}
    +(0.00,2.50) node{}
    +(0.00,0.83) node{}
    +(1.44,0.00) node{}
    +(2.89,0.83) node{}
    +(2.89,2.50) node{}
  };
\end{tikzpicture}
\end{center}
\vskip -0.6 cm \caption{A cubic graph $G$ with $c_{\rm ind}(G)=\frac{n(G)}{2}$.} \label{f:cubic}
\end{figure}

If Conjecture~\ref{conj_cubic} is true, then the bound is achieved, for example, by the graph shown in Figure~\ref{f:cubic}.

\medskip
Although we are only able to prove the NP-completeness for graphs of maximum degree~$4$, we conjecture the following:

\begin{conj}
INDUCED CYCLE is NP-complete for cubic graphs.
\end{conj}

\section{Acknowledgement}

Research of Michael A. Henning supported in part by the South African National Research Foundation and the University of Johannesburg.

\end{document}